\documentclass{article}

\usepackage{amsmath}
\usepackage{amssymb}
\usepackage{amsthm}
\usepackage{amsfonts}
\usepackage{cancel}
\usepackage{mathbbol}
\usepackage{graphicx}
\usepackage[mathscr]{eucal}
\usepackage{txfonts}
\usepackage{bbm,epsfig,graphics,epic,color,rotating,color}

\newtheorem{hypothesis}{Hypothesis}
\newtheorem{theorem}{Theorem}[section]

\newtheorem{definition}[theorem]{Definition}
\newtheorem{remark}[theorem]{Remark}

\newcommand{\reff}[1]{(\ref{#1})}

\renewcommand{\k}{\varkappa}
\newcommand\e{\varepsilon}
\newcommand\la{\lambda}
\renewcommand\phi{\varphi}

\newcommand{\bl}[1]{\mathbf{#1}}
\newcommand{\bs}[1]{\boldsymbol{#1}}

\title{ Hawking-Penrose Black Hole Model. Large Emission Regime. }

\author{ E. Pechersky \thanks{ Supported by FAPESP via grant 2019/08557-0 } 
\\ \small{Institute for Information Transmission Problems, 19, Bolshoj Karetny, Moscow, 127994, RF }
\\ \small{e-mail: pech@iitp.ru } \\ [2ex]
         S. Pirogov \thanks{ Supported by Russian Science Foundation Project 17-11-01098 }
                      \\ \small{ Institute for Information Transmission Problems, 19, Bolshoj Karetny, Moscow, 127994, RF }
                      \\ \small{ e-mail: pirogov@iitp.ru} \\[2ex]
         A. Yambartsev \thanks{ Supported by FAPESP via grant 2017/10555-0 }              
                     \\ \small{Institute of Mathematics and Statistics, University of S\~ao Paulo, S\~ao Paulo, 05508-090,
SP, Brazil} \\ \small{e-mail: yambar@ime.usp.br}
                      }

\begin{document}

\maketitle
\begin{abstract}


In this paper, we propose a stochastic version of the Hawking-Penrose 
black hole model. We describe the dynamics of the stochastic model as a continuous-time Markov jump process of quanta out and in the black hole. The average of the random process 
satisfies the deterministic picture accepted in the physical literature. 
Assuming that the number of quanta is finite the proposed Markov process consists of two
componentes: the number of the quanta in the black hole
and the amount of the quanta outside. 

The stochastic representation allows us to apply large deviation theory to study the asymptotics of probabilities of rare events when the number of quanta grows to infinity. The theory provides explicitly the rate functional for the process. Its infimum over the set of all trajectories leading to large emission event is attained on the most probable trajectory. This trajectory is a solution of a highly nonlinear Hamiltonian system of equations. Under the condition of the stationarity of the fraction of quanta in the black hole, we found the most probable trajectory corresponding to a large emission event.

\end{abstract}


\section{Introduction}

A state of a system is not always the result of a quiet and long evolution. 
Sometimes a very rare event drastically changes  directions of the development. 
If randomness is presented in the system then the rare event can be studied by \textit{large deviations theory}. 
Large deviations theory is one of the well developed and often currently applied parts of the probability theory 
which gives means for asymptotical evaluations of the rare event probability.

Large deviations theory started from the famous Cram\'er's article \cite{Cr}. This work initiated the elaboration of the new section of probability theory. One of the next crucial contributions to the theory has been done by  S. R. S. Varadhan in \cite{Var} where the notion of the \textit{large deviations principle} was introduced. The theory is well developed at present. There exists a fairly large library of books devoted to large deviations theory, \cite{FW, FK, DZ, Rez, Puh, SD, Hol}. 
Some of the books contain chapters about applications of large deviations.

The goal of this paper is an application of large deviation theory to continuous-time Markov processes whose average describes the deterministic evolution of the simplest black hole model considered in \cite{Hak} and \cite{IH}. 


In essence, the concept of a black hole as a domain bounded by an event horizon was discovered by K.Schwarzschild in 1916 \cite{Schw}. The concept had no connection with statistical physics until the calculation of a black hole entropy and the discovery of Hawking radiation, \cite{Bek,Haw1,Haw2}. After this discovery, S. Hawking considered the model containing a black hole and a photon gas in dynamical equilibrium \cite{Hak}. In more detail, this dynamical model was developed by R. Penrose \cite{IH}. 

Here, we propose a stochastic version of the Haw\-king-Penrose model. The study of the stochastic version black hole model is motivated by the quantum nature of the black hole emission \cite{Hak}. Let $N$ be the total number of quanta in the model. The continuous-time Markov process $\xi$ takes its values from 0 to $N$ and describes the evolution of the number of quanta in the black hole (the black hole size). It is formally defined by their infinitesimal generator (\ref{gen}). Since we are interested in the Hawking radiation we introduce into the model the second component, $\eta$, that counts the number of quanta emissions: any hole reduction by one quanta is followed by one emission. The evolution two-component process $\psi=(\xi(\cdot), \eta(\cdot))$ is determined by the generator (\ref{3o1}). The transition rates of this process correspond to the laws of black hole physics.


In this paper we study the large deviations asymptotics for a stochastic version $\psi$ of the Haw\-king-Penrose model 
with special attention to the large emission regime. The emission is a random process that depends on the size of the horizon of the black hole. Because of the randomness there exists a small probability that the average radiation flux during a finite time interval became very large. Despite the small probability, it is strictly positive, thus, such fluctuations will appear with probability one. Stochasticity of the absorption (when the number of quanta in the black hole increases by one) process is due to the non-homogeneity of particle localizations out of the black hole. The physics that causes these fluctuations is not discussed here.

We apply the large deviation theory to find the way how a large emission event occurs. The theory uses an asymptotic approach to the problem. In our case, we study the asymptotics with respect to the growing total number of quanta $N$. Therefore we need to consider the scaled Markov processes $\psi_N=(\xi_N, \eta_N)$, see (\ref{4.8}, \ref{4.6}). The large emission event on the time interval $[0,T]$ is defined in terms of the scaled process: $\eta_N(T) > BT$, where $B$ is the emission rate.

When using the large deviations theory the rate function of the studied stochastic system is sought. To find the rate function of the studied system at appropriate scaling we apply the approach developed in \cite{FK}. The rate function along the asymptotics allows one to find the trajectory of the black hole state which corresponds to the given amount of the emission. When the scaling parameter $N$ is going to infinity the probability is concentrated in every neighborhood of this trajectory. Finding this trajectory is reduced to solving a Hamiltonian system of equations (\ref{6.1}). In the considered case, the Hamiltonian system is highly nonlinear and, unfortunately, finding its solution on the set ${\bl E}$ (see (\ref{4.2})) of all trajectories corresponding to the large emission is a hard problem. We formulate our guess solution as Hypothesis 1.

Taking into account the hypothesis we find the solution on the very restrictive set of trajectories $\bl G$ (see (\ref{G}) for the definition of the set $\bl G$). Namely, we assume that the average hole size is constant and the corresponding emission average is a linear function. We introduce the concept of \textit{stationary emission regime} (see Definition~\ref{Def1}), it is, basically, the solution of Hamiltonian system which belongs to the set $\bl G$.  We prove the following result (see Theorem~\ref{th2}).

\begin{enumerate}
\item[] For each emission rate $B>0$ there exists a mass $m_B$ of the black hole such that the pair of trajectories $(x(t),y(t)) \equiv (x_B, Bt), t\in [0,\mathcal{T}],$ is the stationary emission regime. Here $x_B = \frac{m_B c^2}{E}$, where $E$ is a total energy of the system, and 
$$
m_B \propto \frac{1}{\sqrt[3]{B} },
$$
where the proportionality coefficient is some combination of physical constants.
\end{enumerate}

This relation is new. It describes the correspondence between the size of the black hole and its emission rate in the large emission regime. 

The present work continues our works \cite{PPSVY,PPSVY2,PPSVY3}, where  the similar problems concerning emission regime were studied. The paper is organized as follows. In the next section we recall the deterministic picture, Section~\ref{Deterministic}, and then we formulate our stochastic Markov model in Section~\ref{Stochastic}.  Section~\ref{LD} is devoted to the application of large deviation theory. In this section we derive the rate function (\ref{5.1}). In Section~\ref{Results} we provide the corresponding Hamiltonian system (\ref{6.1}), we formulate the main result, Theorem~\ref{th2}, and the proof. Section~\ref{conclusion} concludes the paper.

\section{Hawking-Penrose black hole model} 

The goal of this paper is to propose and study a stochastic version of the Hawking and Penrose black hole model introduced in \cite{Hak} and \cite{IH}.  The model in our considerations has two constituents: the black hole and a cloud of photons. A part of the photons are located in the hole, the remaining photons are free and located in a box with reflected boundaries. There exists an exchange of photons between the cloud and the hole: emission and absorption. This exchange we describe by a Markov process with discrete phase space.

Before the construction of stochastic model we recall the Hawking-Penrose black hole model accepted in physical literature.

\subsection{Deterministic picture}\label{Deterministic} Let $V$ be a volume with mirror boundaries containing radiation with total energy $E$. Some amount of the energy $e$ is  absorbed by the black hole. The black hole emits a radiation by the Hawking process. It means that the amount of the energy in the black hole depends on time.
\begin{remark}{Remark}
In this subsection we assume that the values of $E, e$ and $m$ take real values. Further, when we will consider the random version, the values of $E, e$ will be discrete.
\end{remark}

The Schwarzschild radius of the black hole equals
\[
R=\frac{2G m}{c^2}=\frac{2G}{c^4}e,
\]
where $m=\frac e{c^2}$ and $G$ is the gravitational constant. The radius $R$ depends on the energy $e$ of the black hole. We denote the coefficient connecting $R$ and $e$ by $a$, 
\begin{equation}\label{RR}
R=ae,
\end{equation}
where
\begin{equation}\label{1.01}
a=\frac{2G}{c^4}.
\end{equation}

The energy $e$  satisfies the balance equation,
\begin{equation}\label{2.4}
\frac{\mathrm{d} e}{\mathrm{d} t}=W_{abs}-W_{em}.
\end{equation}
In this equation the power absorbed by the black hole is 
 \begin{equation}\label{2.5}
 W_{abs}=cA \frac{1}{4} \frac{E-e}{V},
 \end{equation}
  where 
  \begin{equation}\label{2.-}
 A= 4\pi R^2
 \end{equation}
is the horizon area, and the factor 1/4 appears by  geometrical reasons, see \cite{LL}. 
 
 \begin{remark}{Remark} (\textit{i}) The factor 1/4 in \reff{2.5} reflects the fact that the absorbed power falls into black hole at some angle  $\theta$ to the surface. The absorbed power is proportional to $\cos(\theta)$. An average value of $\cos(\theta)$ on the hemisphere $0\leq \theta<\pi/2$ equals to 1/2. Additional factor 1/2 appears because we have to consider only rays directed towards the surface (\cite{Feyn}, Vol.1, Ch. 45). 
 
(\textit{ii}) Here we ignore the gravitational light deflection. An elementary discussion of the gravitational light deflection as a consequence of the equivalence principle discussed in \cite{Ber}, vol.1, Ch. 14.
 
(\textit{iii}) The considering of the light deflection gives $(27\pi/4)R^2$ instead of $\pi R^2$ (see \cite{LL2}, Ch. 12, Section 102: gravitational collapse of spherical body, pp. 338).
\end{remark}

 Using \reff{RR} and \reff{2.-} we obtain
 \begin{equation}\label{2.-1}
 W_{abs}=\frac{\pi c}V a^2e^2(E-e).
 \end{equation}
 The black hole emission $W_{em}$ was calculated in \cite{DW} (eq. (146))
 \[
 W_{em}=\sigma AT^4,
 \]
 where 
 \[T=\frac{\hbar c}{4\pi R}
 \]
 is the Hawking temperature and $\sigma$ is the emission constant (see \cite{DW}, eq.(146)). Note that the emission constant $\sigma$ does not coincide with the classical Stefan-Boltzmann constant \cite{LL}.
Using expressions of $A$ and $R$ via $e$ we obtain
 \[
 W_{em}=\sigma\frac{(\hbar c)^4}{(4\pi)^3a^2}\frac1{e^2}.
 \]
 Let
 \begin{equation}\label{2.1}
 b=\frac{\hbar c}{4\pi a},
 \end{equation}
 then 
 \[
 T=\frac b e.
 \]
 
 We obtain (see \reff{2.4})
 \begin{equation}\label{2.3}
 \frac{\mathrm{d} e}{\mathrm{d} t}=a_1a^2e^2\frac{E-e}{V}-a_2a^2\frac{b^4}{e^2},
 \end{equation}
 where 
 \begin{equation}\label{2.2}
 a_1=\pi c, \ a_2=4\pi\sigma.
 \end{equation}
 
 This equation can have a stationary solution if the equation
 \begin{equation}\label{size}
 a_1e^4(E-e)=a_2b^4V
\end{equation}
 has a solution. The condition for it is
 \begin{equation}\label{size1}
 \frac{4^4}{5^5}E^5\geq \frac{a_2}{a_1}b^4V.
\end{equation}
 If this inequality is strict, then the equation \reff{size} has two solutions. One of them corresponds to the stable and another to the unstable black hole \cite{IH}. 
 
 \subsection{Stochastic picture}\label{Stochastic}  In this section we propose a discrete version of the system  outlined above. Moreover, we impose stochasticity on the system. 
 
 As in the previous section, $E$ is the total energy in the volume $V$, and $e$ is the part of $E$  which is assumed to be contained in the black hole. The discreteness assumes that the total energy $E$  is split in quanta. Let $N$ be the total number of quanta, then the energy $\e$ of each quanta is
 \[
 \e=\frac EN.
 \]
 From now $e$ is also a discrete variable.
 Later on, $E$ is fixed while $N$ is growing. 
 
  The volume $V$ splits into two parts: the black hole interior and its exterior. An arbitrary positive part $k=1,2,...,N$ of the quanta can be absorbed by the black hole, and be contained in it.
  The energy of the black hole is $e=k\e$ if $k$ quanta are in the hole. 
 
 \subsubsection{Markov process} 
The  dynamics consists of  the emission and  the absorption  of the quanta by the black hole. This dynamics we construct as  the continuous-time Markov process $\xi(t),t\in[0,\mathcal{T}]$ with the state space $\mathscr{N}=\left\{1,...,N\right\}$. The state of the process is interpreted as the number of quanta into the black hole.
 The transition rates of $\xi(t)$ are defined as the following:
 
 \begin{enumerate}
\item[] If $\xi(t)=k>1$, then the rate of the transition $k\to k-1$ (the emission rate) equals to $$\frac{W_{em}}{\e}=\frac{a_2a^2b^4}{E^3}N\frac{N^2}{k^2}.$$ 
 
\item[] If $\xi(t)=k<N$, then the rate of the transition $k\to k+1$ (the absorption rate) equals to $$\frac{W_{abs}}{\e}=\frac{a_1a^2E^2}{V}N\frac{k^2}{N^2}\left( 1-\frac kN \right).$$
\end{enumerate}

Thus, the generator of the jump Markov process $\xi(t)$ is
 \begin{equation}\label{gen}
\mathbf{L}f(k)=N\frac{a_1a^2E^2}{V}\frac{k^2}{N^2} \left( 1-\frac kN \right) [f(k+1)-f(k)]+N\frac{a_2a^2b^4}{E^3}\frac{N^2}{k^2} (1-\delta(k-1))[f(k-1)-f(k)].
\end{equation}
Here $\delta(k)=1$ for $k=0$, and $\delta(k)=0$ otherwise. 
\begin{remark}{Remark}
We introduce the term $1-\delta(k-1)$ which does not allow the black hole to evaporate completely.
\end{remark}

We further use the following notations
\begin{equation}\label{7.0}
\begin{aligned} 
\mu &= a_2a^2b^4\tfrac1{E^3}, \\
\la &= a_1a^2\tfrac{E^2}V.
\end{aligned}
\end{equation}

\subsubsection{Markov process with emission} 
Next, we consider the joint process $\psi=(\xi,\eta)$, where the second component  $\eta(t)$ counts the number of quanta emissions from the hole during the time interval $[0, t], \ t\le \mathcal{T}$. The process $\eta(t)$ takes its values in $\mathbb{Z}_+$, and it is non-decreasing process. The initial value $\eta(0)=0$, and we suppose that $\xi(0)$ is uniformly distributed on $\mathscr{N}$. Therefore, the generator of the joint process is
\begin{equation}\label{3o1}
\begin{aligned}
\mathbf{L}f(k,m) &\ = \ \la N\frac{k^2}{N^2}\left(1-\frac kN\right)\left[ f(k+1,m)-f(k,m) \right] \\ 
&\ + \ \mu N\frac{N^2}{k^2} \bigl(1-\delta(k-1) \bigr)\left[ f(k-1,m+1)-f(k,m) \right],
\end{aligned}
\end{equation}
where $k\in\mathscr{N}, m\in\mathbb{Z}_+$.

Considering the large deviations of the black hole emissions during the time interval $[0,\mathcal T]$ we should scale  $(\xi(t),\eta(t))$:
\begin{equation}\label{4.8}
\xi_N(t)=\frac{\xi(t)}N,\ \eta_N(t)=\frac{\eta(t)}N.
\end{equation}
In this scaling we study the large emission when   $N\to\infty$.


The joint process $\psi_N(t)=\left( \xi_N(t),\eta_N(t) \right)$  takes its value in $ D_{N}=\left( \frac{1}{N}\mathscr{N}
\times\frac{1}{N}\mathbb{Z}_+\right) $. Since $ D_{N}\subset D=[0,1] \times\mathbb{R}_+$ for every $ N $ we
will say that the processes $ \psi_{N} $ takes their values in $D$.

The process $\psi_N$ is the jump process with
two types of jumps: $\left(\frac1N,0\right)$ and
$\left(-\frac1N,\frac1N\right)$.  Let $(x_N,y_N)\in D_{N}$. Then the
infinitesimal operator of $\psi_N$ is
\begin{equation} \label{4.6}
\begin{aligned}
\mathbf{L}_{\psi_N}f(x_N,y_N)  &\ = \ \la x_N^2(1-x_N)N \left[ f(x_N+\tfrac1N,y_N)-f(x_N,y_N) \right] \\
&\ + \  \mu\frac1{x_N^2}N \left( 1-\delta(x_N-\tfrac1N)\right) \left[ f(x_N-\tfrac1N,y_N+\tfrac1N)-f(x_N,y_N) \right].
\end{aligned}
\end{equation}

Let $(x,y)\in D$ and a sequence $(x_N,y_N)\in D_N$ be such that $(x_N,y_N)\to (x.y)$.  Assuming differentiability of $f$ we obtain a limit
\[
\mathbf L_\infty f(x,y)=\lim\mathbf{L}_{\psi_N}f(x_N,y_N)=\la x^2(1-x)\frac{\partial f}{\partial x}+\mu\frac1{x^2}\left(\frac{\partial f}{\partial y}-\frac{\partial f}{\partial x}\right).
\]
Let
\begin{equation}\label{4.7}
\bl{S}= \bigl\{(\bl x(t),\bl y(t)):\:[0,\mathcal T]\to  D \bigr\}
\end{equation}
be Skorohod space, it means that the paths $\bl x(\cdot)$ and $\bl y(\cdot)$ are continuous from the right and have limits from the left. This space is equipped with the Skorohod topology \cite{Bil}. Let also $\bl C_1\subset\bl S$ be a subset of pairs of absolutely continuous  functions $(\bl x(\cdot), \bl y(\cdot))$ such that  $\bl x(t)\in[0,1]$ and $\bl y(t)\in\mathbb{R}_+$ is non-decreasing with initial $\bl y(0)=0$.
The process $\psi_N$ induces a measure on $\bl S$.

The operator $\mathbf L_\infty$ can be considered as an infinitesimal generator of a deterministic dynamics, described by the following ordinary differential equations
\begin{eqnarray}
\frac{\mathrm{d} \bl x}{\mathrm{d} t}&=&\la \bl x^2(1-\bl x)-\mu\frac1{\bl x^2} \theta({\bl x}) ,\label{4.4}\\
\frac{\mathrm{d} \bl y}{\mathrm{d} t}&=&\mu\frac1{\bl x^2} \theta({\bl x}), \label{4.5} 
\end{eqnarray}
for $(\bl x,\bl y)\in\bl C_1$, where $\theta(\cdot)$ is the Heaviside step function 
whose value is zero for negative arguments and one for positive arguments, we set $\theta(0)=0$. Note, that the equation \reff{4.4} coincides with the equation \reff{2.3}. The equation \reff{4.5} counts an amount of emitted energy.

For large finite $N$ the paths of a random process $\psi_N$ fluctuate around the solutions of \reff{4.4} and \reff{4.5}.
The probability of these fluctuations are governed by the rate function $I(\bl x(\cdot),\bl y(\cdot))$, which we define in the next section about the large deviation theory. Here we outline the role of the rate function $I$. The probability that the process $\psi_N$ is close to a path  $(\bl x(t),\bl y(t))$ (here $(\bl x(t),\bl y(t))$ does not necessary the solution of \reff{4.4}, \reff{4.5}) has an rough exponential asymptotics
\[
\Pr \Bigl( \psi_N(t)\approx (\bl x(t),\bl y(t)), t\in[0,\mathcal T] \Bigr) \asymp \exp\bigl\{ -NI(\bl x(\cdot),\bl y(\cdot)) \bigr\}
\] 
as $N\to\infty$. The sign $\approx$ means that the process $\psi_N$ is located in a neighborhood of the path $(\bl x(t),\bl y(t))$, and the  neighborhood is shrinking to this path with growing $N$.

\section{Large Deviations}\label{LD}
To find the probability of the large emission on $[0,\mathcal T]$ we use large deviation theory. It is especially useful  when looking at the asymptotic probability  of rare events.
We describe the large deviation approach in terms of the system studied here. The large emission from the black hole on interval $[0,\mathcal T]$ we determine as the event
\begin{equation*}
\mathcal E_N= \bigl( (\xi_N(\cdot),\eta_N(\cdot)) \in \bl S:\ \eta_N(\mathcal T)\geq B\mathcal T \bigr),
\end{equation*}
where $B>0$.  The first component  $\xi_N$ is irrelevant in this event, the same for the values of $\eta_N(t)$ for $t<\mathcal T$ except $t=\mathcal T$, where $\eta_N(\mathcal T)\geq B\mathcal T$. 
Let
\begin{equation}\label{4.2}
{\bl E}=\bigl\{ (\bl x(\cdot),\bl y(\cdot)) \in \bl C_1:\ \bl y(\mathcal T)\geq B\mathcal T, \bl y(0)=0 \bigr\}.
\end{equation}

Further we will consider a more restrictive event $\mathcal G_N \subset \mathcal E_N$ which is related to the so-called \textit{stationary emission regime}, see Definition 1 below. In the definition of $\mathcal G_N$, strong restrictions on the first component $\xi_N(t)$ are introduced, as well. To this end we consider a following subset $\bl G\subset\bl C_1$
\begin{equation}\label{G}
\bl G= \bigcup_{c_1 \in [0,1]} \bigcup_{c_2 \ge B}  \bigl\{ (\bl x(\cdot),\bl y(\cdot))\in \bl C_1:\ \bl x(t)\equiv c_1, \bl y(t)=c_2 t
\bigr\}\subset \bl E.
\end{equation}
Note that both ${\bf E}$ and ${\bf G}$ depend on $B$, but we omit it in notation.

Then, for a given $\delta$ let $U_\delta(\bl G)$ be a $\delta$-neighborhood of $\bl G$ in Skorokhod topology on the space $\bl S$. Finally, the set $\mathcal G_{N,\delta}$ we define as follows:
\[
\mathcal G_{N,\delta}=\bigl( (\xi_N(\cdot),\eta_N(\cdot)) \in U_\delta(\bl G) \bigr).
\]
Thus, the path $(\xi_N(\cdot),\eta_N(\cdot))$ belongs to $\mathcal G_{N,\delta}$ if there exists a trajectory $(\bl x(\cdot),\bl y(\cdot))\in\bl G$ such that $\xi_N(\cdot) \in U_\delta(\bl x(\cdot))$ and $\eta_N(\cdot)\in U_\delta(\bl y(\cdot))$.

\vspace{0.5cm}

The asymptotics  of the probabilities of $\Pr(\mathcal E_N)$ and $\Pr(\mathcal G_{N,\delta})$ as $N\to\infty$ is the subject of the large deviation theory.
The large deviation theory states the existence of the functional
\[
I(\bl x,\bl y):\:\bl C_1 \to \mathbb{R}_+,
\]
such that $I(\bl x,\bl y)=\infty$, when $(\bl x,\bl y)\notin \bl C_1$. In the large deviation theory, the functional $I$ is called the \textit{rate function} which was mentioned in the previous section. The properties of rate function are well described in the literature (see, for example, \cite{FK}).

Applying the large deviations theory \cite{DZ} we can find the logarithmic asymptotics of $\Pr(\mathcal E_N)$ and $\Pr(\mathcal G_{N,\delta})$, that is  
\begin{equation}\label{7.7}
\begin{aligned}
\lim_{N\to\infty}\frac1N\ln \Pr(\mathcal E_N)&=\inf_{(\bl x,\bl y)\in\bl E} I(\bl x,\bl y),\\
\lim_{\delta\to 0}\lim_{N\to\infty}\frac1N\ln \Pr(\mathcal G_{N,\delta})&=\inf_{(\bl x,\bl y)\in\bl G} I(\bl x,\bl y).
\end{aligned}
\end{equation}

Looking for the rate function $I(\bl x,\bl y)$ in our case we follow the method of Feng and Kurtz \cite{FK}. The rate function, according to this method, is constructed by a Hamiltonian $H$. In first step, the non-linear Hamiltonian has to be found: for $(x,y)\in D$ (see \reff{4.6})
\begin{eqnarray*}
({ \cal H}_N f)(x,y)&:=&\frac1N \exp\{-Nf(x,y)\}\times {\mathbf L}_{\psi_N}\exp\{Nf(x,y)\} \\
&=&\la x^2(1-x)\left[\exp\left\{N\left(f\left(x+\tfrac{ 1}N,y\right)-f\left(x,y\right)\right)\right\}-1\right]\\
&& + \ \mu \frac1{x^2}(1-\delta(x-\tfrac1N)\left[\exp\left\{N\left(f\left(x-\tfrac1N,y+\tfrac1N\right)-f\left(x,y\right)\right)\right\}-1\right].
\end{eqnarray*}
It is assumed in the above expression that $0<x<1$ and $N$ is large enough.
Then 
\begin{equation}\label{2.13}
\begin{aligned}
& \lim_{N\to\infty} ({ \mathcal H}_N f)(x,y) \\
& {} =\la x^2(1-x)\left[\exp\{\tfrac{\partial}{\partial x}f(x,y)\}-1\right]+\mu\frac1{x^2}\left[\exp\{-\tfrac{\partial}{\partial x}f(x,y)+\tfrac{\partial}{\partial y}f(x,y)\}-1\right].
\end{aligned}
\end{equation}
Using the notations 
$$
 \k_1:=\tfrac{\partial}{\partial x}f(x,y),\  \k_2:=\tfrac{\partial}{\partial y}f(x,y),
$$
we obtain  from \reff{2.13} the Hamiltonian $H$ of the system
\begin{equation}\label{4_1}
H(x,y,\varkappa_1,\varkappa_2)=\la x^2(1-x)[e^{\varkappa_1}-1]+\mu\frac1{x^2}[e^{-\varkappa_1+\varkappa_2}-1].
\end{equation}

 
 To define the rate function for the considered system we introduce paths $(\bs\k_1,\bs\k_2)$ on $[0,\mathcal T]$; $(\bs\k_1(t),\bs\k_2(t))\in\mathbb{R}^2$.
Then the rate function is obtained as (see  \reff{4_1})
\begin{equation}\label{5.1}
\begin{aligned}
I(\bl x,\bl y) \ = \ \int_0^{\mathcal T} \mathcal{L} (\bl x, \bl y) \mathrm{d} t = \ \int_0^{\mathcal T} \sup_{\bs\varkappa_1(t),\bs\varkappa_2(t)} \Bigl\{ \bs\varkappa_1(t)\dot{\bl x}(t)+\bs\varkappa_2(t)\dot{\bl y}(t)  & \\
 - \la\bl x^2(t)(1-\bl x(t))[e^{\bs\varkappa_1(t)}-1]  
 &  - \mu\tfrac1{\bl x^2(t)}[e^{-\bs\varkappa_1(t)+\bs\varkappa_2(t)}-1] \Bigr\}\mathrm{d} t,
\end{aligned} 
\end{equation}
where
\[
\mathcal{L} (\bl x(t), \bl y(t))=\sup_{\bs\varkappa_1(t),\bs\varkappa_2(t)}\left\{\bs\varkappa_1(t)\dot{\bl x}(t)+\bs\varkappa_2(t)\dot{\bl y}(t)
-H(\bl x(t), \bl y(t),\bs\varkappa_1(t),\bs\varkappa_2(t))\right\}
\]
is Legendre transform of Hamiltonian $H$ \reff{4_1}. Recall that $(\bl x(t), \bl y(t))\in\bl C_1$.

\section{Result}\label{Results}  Our goal is to study how the large emission occurs.
To this end, on the set $\bl E$ of all trajectories that correspond to the large emission we have to find a trajectory where the infimum 
\begin{equation}\label{5.2}
\inf_{(\bl x,\bl y)\in\bl E} I(\bl x,\bl y)
\end{equation}
is attained (see \reff{4.2} for the definition of the set $\bf E$). Note that on the set $\bf E$ there are not any constraints on the fraction of quanta in the black hole. 

Since the rate function $I$ is the non-linear integral functional which integrand is the Legendre transform of Hamiltonian \reff{4_1}, the extremals of \reff{5.2} 
should satisfy a Hamiltonian system
\begin{eqnarray}\label{6.1}
\left\{ \begin{array}{rcl}
\dot {\bl{x}}&=&\la \bl x^2(1-\bl x)\exp\{\bs\varkappa_1\}-\mu\displaystyle\frac1{\bl x^2}\exp\{- \bs\varkappa_1+\bs\varkappa_2\}, \\
\dot {\bl y}&=&\mu\displaystyle\frac1{\bl x^2}\exp\{-\bs\varkappa_1+\bs\varkappa_2\}, \\
\dot{ \bs\varkappa}_1&=&-\la(2\bl x-3\bl x^2) [\exp\{\bs\varkappa_1\}-1]+\mu\displaystyle\frac2{\bl x^3}[\exp\{-\bs\varkappa_1+\bs\varkappa_2\} - 1], \\
\dot{\bs\varkappa}_2&=&0,
\end{array} 
\right.
\end{eqnarray}
with suitable boundary conditions. The system \reff{6.1} is the Euler-Lagrange equation for integral functional $I(\bl x,\bl y)$, see \reff{5.1}.
Due to the high nonlinearity of the system \reff{6.1} we cannot find the solution, but we guess that the minimum is attained on the trajectory which belongs to the set $\bl G$. Thus, the main goal would be the following result which will be formulated as the hypothesis.

\begin{hypothesis}
For any $B>0$ there exists $x_B$ such that the functions 
\[
\bl x_B(t)\equiv x_B\in[0,1],\;\;
 \bl y_B(t)= B t,\,\, t\in[0,\mathcal T]
\]
attain the infimum in \reff{5.2}:
\[
I(\bl x_B,\bl y_B)=\inf_{(\bl x,\bl y)\in\bl E} I(\bl x,\bl y).
\]
\end{hypothesis}

Unfortunately, a proof of this statement is very complicated. But if we restrict the infimum \reff{5.2} on the set $\bl G$, then the proof becomes an easy task, see Theorem~\ref{th2}. The infimum on the restricted set $\bl G \subset \bl E$
\begin{equation}\label{5.2.1}
\inf_{(\bl x,\bl y)\in\bl G} I(\bl x,\bl y),
\end{equation}
also gives the asymptotic behavior of the large emission probability with restrictions on the value of the number of quanta in the black hole, see \reff{G} for the definition of $\bl G$. Namely, in this case the quanta number satisfies periodic boundary conditions on the time interval $[0,\mathcal T]$.
Moreover, in Theorem~\ref{th2} we find the relationship between the size of the black hole and the size of large emission. Before  formulating the next theorem, we introduce the following definition.


\begin{definition}{Definition}\label{Def1} 
For a constant $B>0$, the path $(\bl x_B(t),\bl y_B(t))$ is called a \textsl{stationary emission regime} if 
\begin{enumerate}
\item there is a constant $x_B$ such that $\bl x_B(t)\equiv x_B, \ \ t\in[0,\mathcal T]$,
\item $\bl y_B(t)=B t, \ \ t\in[0,\mathcal T]$,
\item the path $(\bl x_B(t),\bl y_B(t))$ is extremal of $I$  with the boundary conditions $\bl x_B(0)= \bl x_B(\mathcal T)=x_B$ and $\bl y_B(0)=0,\ \bl y_B(\mathcal T)=B\mathcal T$.
\end{enumerate}
\end{definition}

In Theorem~\ref{th2}, we consider a family of all stationary emission regimes (see Definition~\ref{Def1}) which depends on the parameter $B$. The theorem find the relation between $B$ and the value of constant $x_B$ in the stationary emission regime. 


\begin{theorem}{Theorem}\label{th2} For sufficiently large $B>0$, there exists a unique constant $x_B$ which is the root of the equation 
\begin{equation}\label{8.4}
\frac B{1-x_B}-2\mu\frac1{x_B^3}+\la (2x_B-3x^2_B)=0,
\end{equation}
such that the path $\bl x(t)\equiv x_B$,  $\bl y(t)=B t$ is the stationary emission regime. We have $x_B\to 0$ as $B\to \infty$ with the asymptotics 
\[
x_B\sim \frac{\sqrt[3]{2a_2a^2b^4}}{E}\frac1{\sqrt[3]{B}}.
\]
\end{theorem}

\proof We obtain from the definition of the stationary emission regime that 
\begin{eqnarray}\label{7.4}
\left\{ \begin{array}{rcl}
0&=&\la x_B^2(1- x_B)\exp\{\bs\varkappa_1\}-\mu\tfrac1{ x_B^2}\exp\{- \bs\varkappa_1+\bs\varkappa_2\}, \\
B&=&\mu\tfrac1{x_B^2}\exp\{-\bs\varkappa_1+\bs\varkappa_2\}, \\
\dot{ \bs\varkappa}_1&=&-\la(2x_B-3x_B^2) [\exp\{\bs\varkappa_1\}-1]+\mu\tfrac2{ x_B^3}[\exp\{-\bs\varkappa_1+\bs\varkappa_2\} - 1], \\
\dot{\bs\varkappa}_2&=&0.
\end{array} 
\right.
\end{eqnarray}
From the fourth and second equations of \reff{7.4} it follows that $\bs\varkappa_1$ and $\bs\varkappa_2$ do not depend on time. Besides, the following equality 
\begin{equation}\label{7.1}
\la x_B^2(1-x_B)e^{\k_1}=\mu\tfrac1{x_B^2}e^{- \k_1}e^{\k_2}=B,
\end{equation}
where $\k_1\equiv \bs \k_1$ and $\k_2\equiv \bs\k_2$, follows
from the first and second equations of \reff{7.4}. We obtain from these equations
\begin{equation}\label{8.3}
x_B\left[\la x_B(1-x_B)e^{\k_1}-\mu\tfrac1{x_B^3}e^{- \k_1}e^{\k_2}\right]=0.
\end{equation}

Next we prove the equality
\begin{equation}\label{8.2}
\la x_B^2e^{\k_1}-2\mu\frac1{x_B^3}+\la (2x_B-3x^2_B)=0.
\end{equation}
To this end we use the third equation of \reff{7.4} which we rewrite in the following way
\[
-2\la x_B(1-x_B)e^{\k_1}+2\mu\frac1{x_B^3}e^{- \k_1}e^{\k_2}+\la x_B^2e^{\k_1}-2\mu\frac1{x_B^3}+\la (2x_B-3x^2_B)=0.
\]
Using now \reff{8.3}, we obtain \reff{8.2}. Substitute  in \reff{8.2} the value of $e^{\k_1}=\frac B{\la x_B^2(1-x_B)}$ from \reff{7.1} to obtain the equation \reff{8.4}. It is the equation to find $x_B$ via $B$.
Assuming now that $B\to\infty$ we obtain from \reff{8.4} 
\[
x_B\sim \left(\frac{2\mu}{B}\right)^{\tfrac13}
\]
 since $\la$ is a constant and $x_B\in[0,1]$.  \rule{5pt}{5pt}
 
\begin{remark}{Remark} The asymptotics of $x_B$ is determined only by the $\mu$ which depends only on the emission constant $\sigma$ and the coefficient  in the Hawking's formula for the temperature.
\end{remark}


\section{Conclusion}\label{conclusion}

The paper considers a black hole model proposed by Hawking \cite{Hak} and investigated by Penrose \cite{IH}. In addition to the deterministic picture of the black hole dynamics (\cite{Hak}, \cite{IH}),  the random dynamics  driven by a continuous-time Markov process on a finite observation interval $[0,T]$ is introduced. Two characteristics of the black hole are studied in the course of this dynamics: \textit{(i)} the size (volume) of the black hole at every current moment of the observations, and \textit{(ii)} the accumulated value of Hawking emission from the beginning of observations up to the current moment.

The stochasticity permits to consider very rare events that can happen during the stochastic dynamics of the black hole. Here we considered the case when the value of the emission flux far exceeds the average value. The probability of this event was studied from the large deviation point of view. The dynamics of the black hole size at this event is very different from the average. 

The main result was obtained under the additional assumption that the average size of the black hole does not change in the observation time interval under consideration, in which the Hawking emission flux is very large. We proved that the size of the black hole under these assumptions is proportional to $B^{-1/3}$, where $B$ is the total emission over the observation interval $[0, T]$: the greater the total emission $B$, the smaller the hole size.

The construction considered in this work is a Markov random process that describes the stochastic dynamics of a black hole: absorption of matter and Hawking emission. The physics of these phenomena (absorption and emission) is hidden in the stochastic nature of the process. There are many suitable possible Markov processes. We considered a class of Markov processes which average satisfies deterministic behavior in physics. However,  the behavior of the system under rare events might be very sensible to a chosen stochasticity.

\end{document}